\documentclass[10pt]{article}
\usepackage{latexsym,amsfonts,mathrsfs,amsmath}
\protect\baselineskip=15 true pt
\protect
\newtheorem{defi}{Definition}

\newcommand{\Real}{\mathbb R}

 \font\onzemsym=msbm10 scaled \magstep1

\def\Bbb#1{\mbox{{\onzemsym#1}}}

\def\enter{\Bbb Z}

\def\Real{\mathbb{R} }
\begin{document}
%
%

\title{\bf \normalsize ON COMPUTATIONAL ORDER OF CONVERGENCE OF SOME MULTI-PRECISION SOLVERS OF NONLINEAR SYSTEMS OF EQUATIONS}
\author{{\large Miquel Grau-S\'{a}nchez, \`{A}ngela Grau, Jos\'{e} Luis D\'{\i}az-Barrero} \\
                  \hfill    \\
 {\footnotesize {\it Technical University of Catalonia, Departments of Applied Mathematics II and III.}}\\
      {\footnotesize {\it Campus Nord, Jordi Girona 1-3, 08034 Barcelona, Spain. }} \\
      {\footnotesize  E-mail address: miquel.grau@upc.edu, angela.grau@upc.edu, jose.luis.diaz@upc.edu}
      }
\date{}
\maketitle

%
%
%
%

%
\vspace{-7mm}
\begin{abstract}
\noindent In this paper the local order of convergence used in iterative methods to solve nonlinear systems of equations is revisited, where shorter alternative analytic proofs of the order based on developments of multilineal functions are shown. Most important, an adaptive multi-precision arithmetics is used hereof, where in each step the length of the mantissa is defined independently of the knowledge of the root. Furthermore, generalizations of the one dimensional case to $m$-dimensions of three approximations of computational order of convergence are defined. Examples illustrating the previous results are given.

\vspace{4mm}
\end{abstract}

%
%
%
%

{\small
\noindent  \hspace{7mm} {\it Mathematics Subject Classification}: {65H10, 41A25}
}

\section{Introduction}

Well-known analytic techniques  \cite{trau}--\cite{GS} to prove the local order of convergence of iterative methods to solve single nonlinear equations are generalized to systems of nonlinear equations. More precisely, a new proof of the local order for an iterative method to solve the system of equations $F(x)=0$, where $ F:{\Real^m} \longrightarrow {\Real^m}$, is presented. The basic tools used are the formal developments of the function $F$, its inverse and its derivatives in power series. The vectorial expression of the error equation obtained carrying out this procedure, as we will see later on, is
$$
e_{n+1} = H \left( F ' (\alpha),\,F ''(\alpha), \ldots \right) \, e^{\,\rho}_n + \, O \left( e^{\rho +1}_n \right) ,
$$
where $\, \alpha$ is a simple root of $\,F(x)=0$.

\vspace{3mm}
The preceding technique to prove the local convergence is illustrated with several examples in which generalizations of the unidimensional case to $m$-dimensions are performed using the following definitions of computational order of convergence:

\vspace{-2mm}
\begin{itemize}
\item Computational order of convergence (COC) in Weerakoon et al. (2000, \cite{wefe}).

\vspace{-2mm}
\item Approximated order of convergence (ACOC) in Hueso et al. (2009, \cite{HMT2}).

\vspace{-2mm}
\item Extrapolated order of convergence (ECOC) in Grau et al. (2009, \cite{GG2}).

\end{itemize}

\section{Notation and basic results}

To obtain the vectorial equation of the error we need some known results that for ease reference are included in the following. Let $ F:D\subseteq {\Real^m} \longrightarrow {\Real^m}$ be sufficiently differentiable (Fr\'{e}chet-differentiable) in $D$, and therefore with its differentials continuous. If we consider the $k$th derivative of $F$ at $ a \in {\Real^m}$, we have the $k$-linear function

\vspace{-5mm}
\begin{eqnarray*}
F^{(k)} (a) : \Real^m \times \stackrel{{\tiny{k}}}{\breve{\cdots}} \times \Real^m & \longrightarrow  & \Real^m \\
               (h_1, \ldots , h_k) \; \; \; & \longmapsto & F^{(k)} (a) \, (h_1, \ldots , h_k).
\end{eqnarray*}
That is, $F^{(k)} (a) \, (h_1, \ldots , h_k) \in  \Real^m$. It has the following properties:
\begin{enumerate}
\item[P1.]  $\quad F^{(k)} (a) \, (h_1, \ldots , h_{k-1}, \cdot \, ) \in {\mathscr {L}} \left(\Real^m,\,\Real^m\right) \equiv {\mathscr {L}} \left(\Real^m\right)$.
\item[P2.]  $\quad F^{(k)} (a) \, (h_{\sigma(1)}, \ldots , h_{\sigma(k)}) = F^{(k)} (a) \, (h_1, \ldots , h_k)$,
    where $\sigma$ is any permutation of the set $\{1, \,2, \ldots k\}$.
\end{enumerate}

Notice that from P1 and P2 we can use the following notation:
\begin{enumerate}
\item[N1.]  $\quad F^{(k)} (a) \, (h_1, \ldots , h_{k}) =  F^{(k)} (a) \, h_1  \cdots  h_{k} \, $.
\item[N2.]  $\quad F^{(k)} (a) \, h^{k-1} \; F^{(\ell)} (a) \, h^{\ell} = F^{(k)} (a) \,F^{(\ell)} (a)\: h^{k + \ell -1} \, $.
\end{enumerate}

Hence, we also can express $ F^{(k)} (a) \, (h_1, \ldots , h_{k})$ as
$$
F^{(k)} (a) \, (h_1, \ldots , h_{k-1})\, h_k = F^{(k)} (a) \, (h_1, \ldots , h_{k-2}) \; h_{k-1}\,h_k = \ldots = F^{(k)} (a) \, h_1  \cdots  h_{k} \, .
$$

On the other hand, for any $\,q = \alpha + h \in {\Real^m}$ lying in a neighborhood of a simple zero,
$\alpha \in {\Real^m}$, of the system $\, F (x) = 0$ we can apply Taylor's formulae and assuming that there exists
$\Gamma = \left[F'(\alpha)\right]^{-1}$, we have

\vspace{-2mm}
\begin{equation}\label{eq01}
    F(\alpha + h) = F'(\alpha) \left[ \, h +\, \sum_{k=2}^{3} \,A_k\,h^k +\,O_4 \,\right] ,
\end{equation}

\vspace{-3mm}
where
$$
 A_k = \frac{1}{k!} \; \Gamma \: F^{(k)}(\alpha) , \; \; k \ge 2 , \; \; {\mbox{and}} \; \;  O_4 = O( h^4) .
$$

\vspace{1mm}
Note that since $\, F^{(k)}(\alpha) \in \mathscr {L} \left(\Real^m \times \stackrel{{\tiny{k}}}{\breve{\cdots}} \times \Real^m,\,\Real^m\right)$ and $\,  \Gamma  \in \mathscr {L}\left(\Real^m \right)$, then
$ A_k\, h^k \in \Real^m$.
Moreover, we can express the differential of first order as:

\vspace{-2mm}
\begin{equation}\label{eq02}
    F'(\alpha + h) = F'(\alpha) \left[ \,I +\, \sum_{k=2}^{3} \,M_k\, h^{k-1} +\,O_3 \,\right] ,
\end{equation}
 where $I$ is the identity and $\, M_k = k\,A_k \,$. Therefore, $M_k\,h^{k-1} \in  {\mathscr {L}} \left(\Real^m\right)$. From  (\ref{eq02}), we get
\begin{equation}\label{eq03}
\left[ F'(\alpha + h) \right]^{-1} =  \left[ I - Q_2\,h + Q_3\,h^2 + O_3 \right] \, \Gamma ,
\end{equation}

\vspace{-3mm}
where
$$
\left\{
\begin{tabular}{rcl}
 $Q_2$ &=& $M_2$, \\ [0.6ex]
 $Q_3$ &=& $M^2_2 - M_3$. \\
\end{tabular}
\right.
$$
Hereof we have used the following notation:  $\: M_{k\ell}\,h^{k+\ell-2} = M_k\,h^{k-1}\:M_{\ell}\,h^{\ell-1}$, if $k \ne \ell$, and  $\: M^2_{\ell}$ instead of $M_{\ell\ell}$, if $k = \ell$ .

\vspace{2mm}
In general, if
$
\left[ F ' (d) \right]^{-1} = \left[ I - Q_2 \, \delta + Q_3\,  \delta^2 \right] \, \Gamma \:
$ and
$
 \; F  (t)  =  F ' (\alpha)  \, \left[ \tau + A_2 \, \tau^2 + A_3 \, \tau^3 \right] \,
$,
with the same notation, we obtain:
\begin{eqnarray} \label{eq04}
\left[ F ' (d) \right]^{-1} \: F  (t)  &=&  \tau + \left( A_2 \, \tau^2 - Q_2 \,  \delta \, \tau \right) +
            \left( A_3 \,  \tau^3 - Q_2 A_2 \,  \delta\,  \tau^2 + Q_3\,  \delta^2\,  \tau \right) .
\end{eqnarray}



\vspace{1mm}
We close this section applying the preceding to the Newton's method. Let $x \in {\Real^m}$, Newton's method is given by
\begin{equation}\label{New}
 z  =   x - \, [F'(x)]^{-1} \, F(x).
\end{equation}

Putting $\, q = x $ in (\ref{eq01}) and (\ref{eq03}) we obtain $ F(x)$ and $ [F'(x)]^{-1} $ in powers of $ \, e = x - \alpha$.
The expression of $ \, E = z - \alpha$ in terms of $ e$ is build up subtracting $\alpha$ to both sides of (\ref{New}). Namely,
\begin{equation}\label{Err}
 E = \, e -  \, [F'(x)]^{-1} \, F(x) = \, T_2\, e^2 + T_3\,e^3  + O_4 ,
\end{equation}
where $\; T_2 = M_2/2 = A_2$, and $ \; T_3 =  ( 4\,M_3 - 3\,M^2_2 )/6 = 2\,( A_3 - A_2^2 )$. These values agree with the classical asymptotical constant in the one dimensional case.

\vspace{3mm}
Without to use norms we can define the local order of convergence for one-step iterative method as follows. The local order of convergence is  $p \in \mathbb{N}$ if there exists a $p$--linear function
 $ K \in \mathscr {L} \left(\Real^m \times \stackrel{{\tiny{p}}}{\breve{\cdots}} \times \Real^m, \Real^m \right)$ such that
\begin{equation}\label{Ord}
   E = K \: e^p + O_{p+1} \, ,
\end{equation}
where $\,e^p  $ is $(e, \stackrel{{\tiny{p}}}{\breve{\cdots}}, e)$ and $\; O_{p+1} = O(e^{p+1})$.


\section{Variants of Newton's method}

In the following, using the results presented in the previous section,
three known variants of Newton's method with local order of convergence equal three, are analyzed. We explicitly give their vectorial error equation in which it appears a $3$--linear application instead of the asymptotical error constant used in the one dimensional case.

\subsection{Arithmetic mean Newton's method}
The first variant (\cite{trau},\cite{wefe}), that substitutes the derivative of $F(x)$ by the arithmetic mean of the derivatives of $F$ at the points $x$ and $z$ is given by
\begin{equation}\label{AMN1}
 X = \,x -\, 2 \; [F '(x) + F '(z)]^{-1} \: F(x).
\end{equation}

\vspace{2mm}
From (\ref{Err}) and the developments of $F ' (x)$, $F ' (z)$, we derive the development of $ \left[ F ' (x) + F ' (z) \right]^{-1}$. Indeed,

\vspace{-5mm}
\begin{eqnarray*}
F ' (x) & = & F ' (\alpha) \, \left[ I + M_2 \,e + M_3\, e^2 \right] + O_3 ,\\ [0.8ex]
F ' (z) & = & F ' (\alpha) \, \left[ I + M_2\, E  \right] + O(E^2)\\  [0.6ex]
        & = & F ' (\alpha) \, \left[ I + M_2 \, T_2\, e^2 \right] +O_3 ,
\end{eqnarray*}

\vspace{-1mm}
where $\; T_2 = M_2/2$. From (\ref{eq03}) we get
\begin{eqnarray}\nonumber
 \left[ F ' (x) + F ' (z) \right]^{-1} &=&  \left[ F ' (\alpha) \, \left\{ 2\,I + M_2 \,e +\,\left( M_3 + M_2 \, T_2 \right)  e^2 \right\} \right]^{-1} \\ \label{AMN2}
    &=&
 \frac{1}{2} \left[ I - \frac{1}{2} \left( M_2\, e +  M_3\,  e^2 \right) \right] \,  \Gamma \,+ O_3 .
\end{eqnarray}

Taking into account that $\,M_k = k\,A_k , \; \; k=2,3\;$, subtracting $\alpha$ from (\ref{AMN1}), and  applying  (\ref{eq04}) we have
\begin{eqnarray} \nonumber
X - \alpha &=&\left( - A_{3} + \frac{1}{2}\,{M_{{2}}} \, A_2 +  \frac{1}{2}\,M_{3} \right) e^3 + O_4 \\\label{AMN3}
          &=& \left( \frac{1}{2}\,A_{{3}}+ A_{2}^{2} \right) {e}^{3} + O_4 .
\end{eqnarray}

Note that we get the same expression as the one obtained for the one dimensional case in \cite{trau}, \cite{wefe} and  \cite{GD5}.


\subsection{Harmonic mean Newton's method}
The second variant (\cite{trau},\cite{home},\cite{ozba}), that substitutes the derivative of $F(x)$ by the harmonic mean of the derivatives of $F$ at the points $x$ and $z$ is given by

\begin{equation}\label{HMN1}
 X = \,x -\, \frac{1}{2} \: \left[ F'(x)^{-1} + F'(z)^{-1} \right] \, F(x).
\end{equation}

\vspace{2mm}
From the developments of $F ' (x)$, $F ' (z)$ we derive the developments of $ F ' (x)^{-1}$ and $ F ' (z)^{-1}$. Namely,

\vspace{-6mm}
\begin{eqnarray} \nonumber
F ' (x)^{-1} &=& \left[ I - M_2 e + \left( M^2_2  - M_3 \right) e^2 \right] \,\Gamma + O_3 ,\\  \nonumber
F ' (z)^{-1} &=& \left[ I - M_2 E  \right]\,\Gamma + O(e^4) \\\label{HMN2}
             &=& \left[ I - \frac{1}{2} \: M^2_2 \,e^2 \right]\,\Gamma + O_3 .
\end{eqnarray}

\vspace{-2mm}
From (\ref{HMN2}) we have
\begin{equation}\label{HMN3}
\frac{1}{2} \left[ F'(x)^{-1} + F'(z)^{-1} \right] = \left[ I - \frac{1}{2}  M_2\,e + \frac{1}{2} \: \left( \frac{1}{2} M^2_2 - M_3 \right) e^2  \right] \,\Gamma + O_3 .
\end{equation}

Finally, subtracting  $\alpha$ from (\ref{HMN1}),taking into account (\ref{HMN3}) and applying (\ref{eq04}), yields

\vspace{-4mm}
\begin{eqnarray} \nonumber
    X - \alpha &=& \left( - A_3 + \frac{1}{2}\: M_2 \,A_2 - \frac{1}{2}\: \left( \frac{1}{2}\:M^2_2 -
    M_3 \right)\right) e^3 + O_4 \\ \label{HMN4}
     &=&  \frac{1}{2}\: A_3 \, e^3 + O_4 .
\end{eqnarray}

This result agrees with the ones obtained in \cite{trau}, \cite{home} and \cite{ozba}.


\subsection{Frozen derivatives in Newton's method}

The third iterative method presented in this section is derived independently in the works of Shamanskii (\cite{sham}, 1967) and Potra et al. (\cite{popt}, 1989). It is defined by
\begin{equation}\label{FDN1}
 X = \,z -\: F'(x)^{-1} \, F(z) ,
\end{equation}
where $z = x - \, F'(x)^{-1}\, F(x) $ is a Newton point. This method is a modification of Newton's iterative function where $F $ is computed in the second step without evaluating $F '$. So, we can consider that the derivative is frozen in this second step (\cite{trau},\cite{orrh}--\cite{ABG}). Note that only one computation of the inverse function is needed and it is only necessary one LU decomposition.

\vspace{2mm}
Subtracting $\alpha$ from (\ref{FDN1}) and applying (\ref{eq04}) we obtain the following vectorial error equation:
\begin{eqnarray} \nonumber
X - \alpha &=& E - F'(x)^{-1} \, F(z) \\ \nonumber
           &=& T_2 \,e^2 + T_3\,e^3 - \left[ I- M_2\, e + O_2 \right] \,
           \left[ A^2_{2} \,e^2 + 2 \left( A_3 -A^2_2 \right) e^3 + O_4 \right] + O_4 \\ \label{FDN2}
          &=& 2\, A_{2}^{2} \, {e}^{3} + O_4 .
\end{eqnarray}


\section{Main numerical results}

In this section we introduce two variants of the definition of COC which are independent of the knowledge of the root and we give some ways to compute with adaptive multiprecision arithmetics and a stopping criteria. Moreover, we apply all the results obtained yet to several examples.

\subsection{Theoretical concepts}

Let  $x_{n-2}, \, x_{n-1}, \, x_{n}$ and $x_{n+1}$ be the last four consecutive iterations of the sequence $\{x_n\}_{n\ge 0}$, where $x_n \in \Real$. The definitions of Computational
Order of Convergence (COC),$\; {\bar{\rho}_n}$ \cite{wefe}, Approximated Computational Order of Convergence (ACOC),$\; {\hat{\rho}_n}$ \cite{HMT2}, and Extrapolated Computational Order of Convergence (ECOC),$\; {\tilde{\rho}_n}$ \cite{GG2}, are the following

\vspace{1mm}
\begin{defi} The COC, ACOC and ECOC of a sequence $\{x_n\}_{n\ge 0}$
is defined by

$$
\vspace{-2mm}
\begin{array}{ccc} \label{COC}
     {\bar{\rho}_n} = {\displaystyle \frac{\ln \left| e_{n+1}/e_n \right|}{\ln \left| e_n/e_{n-1} \right|},} & \quad
      \hat{\rho}_n = {\displaystyle
 \frac{\ln \,\left|\, \hat{ e}_{n+1} / \hat{ e}_{n} \right|}{\ln \, \left|\, \hat{ e}_{n} /
 \hat{ e}_{n-1} \right|} \,} ,& \quad
   {\tilde \rho_n} = {\displaystyle
 \frac{\ln \,\left| \,{\tilde e}_{n+1} / {\tilde e}_{n} \right|}{\ln \, \left| \,{\tilde e}_{n} /
 {\tilde e}_{n-1} \right| } \,}
\end{array}
$$

\vspace{2mm}
\noindent respectively, where $e_n =  x_n - \alpha$, $\, \hat{ e}_n =   x_n - x_{n-1} \,$, ${\tilde e}_n =  x_n - {\tilde \alpha}_n  \,$ and $\,{\tilde \alpha}_n = x_n - \, \frac{ \left(\Delta\, x_{n-1}\right)^2}{\Delta^2 \, x_{n-2}} , \; \: n \ge 2$, where $\Delta$ is the forward difference operator, $\, \Delta\, x_k \,= \, x_{k+1} - x_k$.
\end{defi}

\vspace{2mm}
\noindent  One of the main drawback of the COC is that it involves the exact root $\alpha$, which in a real situation it is not known a priori. To avoid this, we have introduced these two variants of COC, that do not use the exact root.

\vspace{2mm}
\noindent In numerical problems where a huge number of significant digits of the solution is needed it
is required the use of methods with a high order of convergence together with an adaptive arithmetics that is to update the length of the mantissa at each step by means of the formula

\vspace{-1mm}
\begin{equation}\label{eqdig1}
    {\tt Digits :=} \; \left[ \, \rho \times \left( - \log | \,e_n | +\, j \right) \right] \, ,
\end{equation}
where $\rho$ is the order of convergence of the method and $[ x]$ denotes the integer part of  $ x $. Notice that the length of the mantissa is increased approximately by the order of convergence $\rho$.
We have numerically checked the value of $j$, by varying it between  $1$ and $5$, in order to have
enough accuracy in the computation of the iterates  $\{x_n\}_{n\ge0}$. We have realized that the minimum value that guarantees all the significant  digits required is, in almost cases,  $j = 2$. Consequently, hereof we consider  $ j = 2 $ in formula (\ref{eqdig1}) except if another value is explicitly given. In addition, to compute $e_n$, $\hat{ e}_n$ or $\tilde{ e}_n$ with an appropriate number of figures, using Definition 1, we must to enlarge the mantissa in the computation of
$\, x_{n+1}, x_n, x_{n-1},\dots $ with at least four additional significant digits.

\vspace{1mm}
\noindent Two relationships, one between $\, e_n$  and $\hat{e}_n$, and the other between $\, e_n$ and  $\tilde{e}_n$ are given in \cite{GNG}. Namely,

\vspace{-6mm}
\begin{eqnarray} \label{errors}
e_n \approx \: C^{\:\frac{1}{1-\rho}} \, \left( \frac{{\hat e}_{n}}{{\hat e}_{n-1}} \right)^{\rho^2 / (\rho -1)} & {\rm and {\hspace{2mm}} } & e_n \approx C^{\:\frac{\rho-1}{2\rho-1}} \:\;
{\tilde e}_n ^{\; \rho^2 /\, (2\rho-1)} .
\end{eqnarray}


\vspace{1mm}
\noindent Notice that for updating the adaptive arithmetic process (\ref{eqdig1}) it is necessary to
know the exact root $\alpha $. In this case the following stopping criteria is applied:

\vspace{-2mm}
\begin{equation}\label{stop1}
     | e_n | = |x_n - \alpha|< 0.5 \cdot 10^{-\eta},
\end{equation}
where $\eta$ is the number of correct decimals and $0.5 \cdot 10^{-\eta}$ is the required accuracy.
The result given in  (\ref{errors})  allows us to substitute the error in (\ref{eqdig1}) by an expression that does not involve the exact root. Indeed, we implement the following adaptive multi-precision arithmetic schemes:

\vspace{-3mm}
\begin{eqnarray} \label{eqdig3}
\hspace{-7mm} {\tt Digits } := \left[ \frac{\rho^3}{\rho-1} \times \left(- \log \left| \, \delta_n \right| +\, j \right) \right] & {\rm or} \: & {\tt Digits } := \left[ \frac{\rho^3}{2\rho-1} \times \left(- \log | \,\tilde e_n| +\, j \right) \right],
\end{eqnarray}
where $\, {\displaystyle \delta_n =  \left| \, \frac{{\hat e}_{n}}{{\hat e}_{n-1}} \right| }$. Moreover, from (\ref{errors}) we propose the following stopping criteria, instead of (\ref{stop1}):

\vspace{-3mm}
\begin{eqnarray}\label{stop2}
    \delta_n < 0.5 \cdot 10^{-\eta\: (\rho -1) / \rho^2}
     & {\rm or} \; &
      | \, {\tilde e}_n  | < 0.5 \cdot 10^{-\eta\: (2 \rho -1) / \rho^2}.
\end{eqnarray}

%
%

\subsection{Examples}

We generalize the preceding definition and techniques for solving seven systems of nonlinear equations using the Maple computer algebra system. We use the norm  $\| \cdot \|_{\infty}$ instead of the absolute value in (\ref{eqdig1})--(\ref{stop2}).

\vspace{2mm}
\noindent We have computed the solution of each system for the same set of initial approximations $\, x^i_0, \: i = 1,2,3$, which have been chosen close to the root $\alpha$ using the eucliden distance $ d_i = \| x^i_0 - \alpha \|_2 $, and taking into account the value of $ \,D_i = \| F(x^i_0) \|_{\infty}$.

\vspace{2mm}
\noindent Depending on the computational order of convergence used, COC, ACOC  or ECOC, the iterative method was stopped when the condition (\ref{stop1}) or (\ref{stop2}) is fulfilled. Note that
in all cases  $\, \eta = 2800$ and we also obtain $\, \| F(x_k) \| < 0.5 \cdot 10^{-\eta}$. Finally, we choose $j$,     see (\ref{eqdig1}) and (\ref{eqdig3}), such that $\, \|e_k\|$, $ \|{\hat{e}}_k\|$ or $\, \| {\tilde{e}}_k\|$ respectively, and $\, \| F(x_k) \|$ have three significant digits al least.

%
%
%
%

\vspace{2mm}
\noindent Tables 1--7 show, for each method and each function, the number, $k$, of iterations needed to compute the root to the level of precision described. Note that independently of using (\ref{eqdig1}) or (\ref{eqdig3})  the number of necessary iterations to get the desired precision is the same. In addition, in the sixth, eighth and tenth column, it is shown an error bound for the corresponding Computational Orders of Convergence (COC, ACOC and ECOC), given respectively by $\: \bar{\rho}_{k-1} = \rho \pm \Delta \bar{\rho}_{k-1}$, $\: \hat{\rho}_k =  \rho \pm \Delta \hat{\rho}_k \,$, and $\: \tilde{\rho}_k =  \rho \pm \Delta \tilde{\rho}_k \,$.

%
%
%
%

\subsubsection{Example 1}

We begin with the system $F_1(x) = 0$ defined by

\begin{center}
$
 \displaystyle \left\{\begin{array}{r} e^x-2=0,\\
\sin(2y-x)=0.\end{array}\right.\
$
\end{center}

\vspace{1mm}
The roots of $F_1(x)=0$ are $(x,y)=(\ln{2},k\pi+\ln{\sqrt{\mathstrut 2}})^t\,$\ with $k\in\enter$.
We study the convergence of iterative methods previously presented towards the root
$\, \displaystyle \alpha ={(\ln{2},\, \ln{\sqrt{2}})}^t \approx {(0.6931471806, \, 0.3465735903)}^t\,.$

\vspace{1mm}
The initial approximations of the four methods are
$x_0^1 =(1,\, 0)^t$ with $d_1 = 0.347$ and $D_1=0.841$\,; \
$x_0^2 =(0.6,\, 0.3)^t$ with $d_2=0.0931$ and $D_2=0.178$\,; \  and \
$x_0^3 =(0.7,\, 0.35)^t$ with $d_3=0.00685$ and $D_3=0.0137$\,.

\begin{table}[ht]
\caption{Numerical results for system $F_1(x)=0$\,.}
\label{angtable1}
{\scriptsize
\noindent \begin{tabular}{cccccccccc}
         \hline\hline \\[-0.7em]
         &       &          &                      & \hfill COC & & \hfill ECOC & & \hfill ACOC & \\[0.2em]
 Method  & $x_0$   & $k$ & ${\|F_1(x_k)\|}_{\infty}$ & ${\|e_{k-1}\|}_{\infty}$             & $ \Delta{\overline{\rho}}_{k-1} $
                                                   & ${\|{\widetilde{e}}_{k}\|}_{\infty}$ & $ \Delta{\widetilde{\rho}}_k $
                                                   & $\delta_{k}$                         & $ \Delta{\widehat{\rho}}_k $
                                                    \\[0.2em]  \hline \\
    NM   & $x_0^1$  & $12$ & $1.02\cdot 10^{-3429}$ & $3.19\cdot 10^{-1715}$ & $4.03\cdot 10^{-15}$ & $1.19\cdot 10^{-2572}$
                           & $3.78\cdot 10^{-19}$   & $3.99\cdot 10^{-858}$  & $6.74\cdot 10^{-20}$ \\[0.6em]
$\rho=2$ & $x_0^2$  & $12$ & $2.95\cdot 10^{-5428}$ & $1.72\cdot 10^{-2714}$ & $5.17\cdot 10^{-16}$ & $1.59\cdot 10^{-4071}$
                           & $4.16\cdot 10^{-20}$   & $9.27\cdot 10^{-1358}$ & $1.44\cdot 10^{-25}$ \\[0.6em]
$(\ast)$ & $x_0^3$  & $11$ & $8.63\cdot 10^{-5050}$ & $2.94\cdot 10^{-2525}$ & $1.91\cdot 10^{-15}$ & $1.13\cdot 10^{-3787}$
                           & $7.61\cdot 10^{-19}$   & $3.83\cdot 10^{-1263}$ & $3.32\cdot 10^{-25}$ \\[0.2em]
\hline \\
   AMN   & $x_0^1$  & $9$  & $1.04\cdot 10^{-7959}$ & $4.93\cdot 10^{-2832}$ & $2.69\cdot 10^{-17}$ & $1.91\cdot 10^{-4422}$
                           & $1.40\cdot 10^{-13}$   & $1.65\cdot 10^{-1769}$ & $2.34\cdot 10^{-35}$ \\[0.6em]
$\rho=3$ & $x_0^2$  & $8$  & $3.46\cdot 10^{-8274}$ & $1.73\cdot 10^{-2758}$ & $9.86\cdot 10^{-18}$ & $3.73\cdot 10^{-4597}$
                           & $1.83\cdot 10^{-17}$   & $2.15\cdot 10^{-1839}$ & $9.74\cdot 10^{-36}$\\[0.6em]
         & $x_0^3$  & $7$  & $3.13\cdot 10^{-5256}$ & $1.68\cdot 10^{-1752}$ & $9.86\cdot 10^{-18}$ & $1.64\cdot 10^{-2920}$
                           & $4.73\cdot 10^{-18}$   & $9.78\cdot 10^{-1169}$ & $2.01\cdot 10^{-35}$\\[0.2em]
\hline\\
   HMN   & $x_0^1$  & $9$  & $6.80\cdot 10^{-7162}$ & $1.01\cdot 10^{-2387}$ & $5.80\cdot 10^{-20}$& $3.25\cdot 10^{-3979}$
                           & $2.05\cdot 10^{-19}$   & $3.23\cdot 10^{-1592}$ & $5.00\cdot 10^{-35}$\\ [0.6em]
$\rho=3$ & $x_0^2$  & $7$  & $2.83\cdot 10^{-3451}$ & $1.19\cdot 10^{-1150}$ & $3.26\cdot 10^{-20}$ & $1.26\cdot10^{-1917}$
                           & $2.63\cdot 10^{-21}$   & $1.06\cdot 10^{-767}$  & $5.35\cdot 10^{-35}$\\ [0.6em]
         & $x_0^3$  & $7$  & $7.52\cdot 10^{-5912}$ & $7.67\cdot 10^{-1971}$ & $5.00\cdot 10^{-18}$ & $1.30\cdot 10^{-3284}$
                           & $3.06\cdot 10^{-21}$   & $1.70\cdot 10^{-1314}$ & $1.72\cdot 10^{-31}$\\ [0.2em]
\hline\\
   FDN   & $x_0^1$  & $8$  & $2.52\cdot 10^{-4653}$ & $1.36\cdot 10^{-1551}$ & $3.12\cdot 10^{-32}$
                           & $2.98\cdot 10^{-862}$ & $5.02\cdot 10^{-16}$ & $9.75\cdot 10^{-1035}$ & $2.52\cdot 10^{-117}$\\ [0.6em]
$\rho=3$ & $x_0^2$  & $8$  & $5.46\cdot 10^{-7653}$ & $1.76\cdot 10^{-2551}$ & $2.42\cdot 10^{-17}$ & $4.39\cdot10^{-4252}$
                           & $3.10\cdot 10^{-18}$   & $2.49\cdot 10^{-1701}$ & $9.24\cdot 10^{-36}$\\ [0.6em]
$(\star)$ & $x_0^3$ & $7$  & $9.43\cdot 10^{-5065}$ & $9.81\cdot 10^{-1689}$ & $9.63\cdot 10^{-18}$ & $3.57\cdot 10^{-2814}$
                           & $1.35\cdot 10^{-17}$ & $3.64\cdot 10^{-1126}$ & $4.06\cdot 10^{-35}$\\ [0.2em]
\hline\hline
\end{tabular}}
\end{table}

\vspace{2mm}
Note that in the results shown in Table \ref{angtable1} in the case marked with $(\ast)$ it is necessary to take $j=4$  for $x_0^1$ and $x_0^2$ in (\ref{eqdig1}) and (\ref{eqdig3}) because for $j<4$ we loose some significant figures in the computation of $F_1(x_k)$. In the case marked with $(\star)$, for the same reason, we take $j=8$ for $x_0^3$.

\subsubsection{Example 2}

The second example involves two quadratic polynomials. Namely,  $F_2(x)=0 $, defined by

\begin{center}
$\displaystyle\left\{
\begin{array}{r} x^2-4x+y^2=0,\\ 2x+y^2-2=0.\end{array}\right.\ $
\end{center}

\vspace{1mm}
Its solutions are  $(x,y)= \left(3-\sqrt{\mathstrut 7}, \pm\sqrt{\mathstrut -4+2\sqrt{\mathstrut 7}}\right)^t$ and
$(x,y)= \left(3+\sqrt{\mathstrut 7}, \pm\sqrt{\mathstrut -4-2\sqrt{\mathstrut 7}}\right)^t$.
We test the convergence of the methods towards the root
$\displaystyle \alpha = \left( 3-\sqrt{\mathstrut 7}, \sqrt{\mathstrut -4+2\sqrt{\mathstrut 7}} \right)^t \approx(0.3542486889, 1.136442969)^t\,.$

\vspace{2mm}
\noindent The initial values are
$x_0^1 =(-1,\,0.4)^t$ with $d_1=1.354$ and $D_1=5.16$\,; \
$x_0^2 =( 0,\,1)^t$  with $d_2=0.354$ and $D_2=1.0$\,; \ and
$x_0^3 =(0.3,\,1.1)^t$ with $d_3=0.0542$ and $D_3=0.19$\,.

\begin{table}[ht]
\caption{Numerical results for system $F_2(x)=0$\,.}
\label{angtable2}
{\scriptsize \noindent
\begin{tabular}{cccccccccc}
         \hline\hline \\[-0.7em]
         &       &          &                      & \hfill COC & & \hfill ECOC & & \hfill ACOC & \\[0.2em]
    Method  & $x_0$   & $k$ & ${\|F_2(x_k)\|}_{\infty}$ & ${\|e_{k-1}\|}_{\infty}$             & $ \Delta{\overline{\rho}}_{k-1} $
                                                      & ${\|{\widetilde{e}}_{k}\|}_{\infty}$ & $ \Delta{\widetilde{\rho}}_k $
                                                      & $\delta_{k}$                         & $ \Delta{\widehat{\rho}}_k $
                                                        \\[0.2em] \hline \\
 NM   & $x_0^1$  & $14$  & $4.33\cdot 10^{-3666}$ & $2.08\cdot 10^{-1833}$ & $1.68\cdot 10^{-14}$ & $3.78\cdot 10^{-1375}$
                               & $2.35\cdot 10^{-15}$   & $3.03\cdot 10^{-917}$  & $4.98\cdot 10^{-40}$ \\[0.6em]
 $\rho=2$ & $x_0^2$  & $12$  & $3.34\cdot 10^{-3778}$ & $1.83\cdot 10^{-1889}$ & $2.38\cdot 10^{-14}$ & $3.43\cdot 10^{-1417}$
                               & $2.98\cdot 10^{-15}$   & $5.33\cdot 10^{-473}$  & $3.26\cdot 10^{-40}$ \\[0.6em]
  $(\ast)$  & $x_0^3$  & $11$  & $3.46\cdot 10^{-3375}$ & $5.88\cdot 10^{-1688}$ & $1.88\cdot 10^{-13}$ & $9.46\cdot 10^{-2532}$
                               & $5.02\cdot 10^{-12}$   & $1.27\cdot 10^{-422}$  & $3.16\cdot 10^{-32}$ \\[0.2em]
\hline \\
   AMN   & $x_0^1$  & $9$  & $6.91\cdot 10^{-5380}$ & $1.16\cdot 10^{-1793}$ & $8.54\cdot 10^{-11}$ & $3.45\cdot 10^{-2989}$
                           & $3.96\cdot 10^{-18}$   & $2.97\cdot 10^{-1196}$ & $6.45\cdot 10^{-36}$ \\[0.6em]
$\rho=3$ & $x_0^2$  & $8$  & $2.40\cdot 10^{-7591}$ & $8.17\cdot 10^{-2531}$ & $8.68\cdot 10^{-13}$ & $8.90\cdot 10^{-4218}$
                           & $4.23\cdot 10^{-11}$   & $1.09\cdot 10^{-1687}$ & $1.40\cdot 10^{-25}$ \\[0.6em]
         & $x_0^3$  & $7$  & $3.95\cdot 10^{-3817}$ & $9.64\cdot 10^{-1273}$ & $2.77\cdot 10^{-13}$ & $5.45\cdot 10^{-2121}$
                           & $7.54\cdot 10^{-17}$   & $5.65\cdot 10^{-849}$  & $5.59\cdot 10^{-35}$ \\[0.2em]
\hline\\
   HMN   & $x_0^1$  & $7$  & $4.33\cdot 10^{-3666}$ & $5.20\cdot 10^{-1230}$  & $1.06\cdot 10^{-16}$ & $2.81\cdot 10^{-1604}$
                           & $2.76\cdot 10^{-17}$   & $4.08\cdot 10^{-688}$  & $1.20\cdot 10^{-28}$ \\ [0.6em]
$\rho=4$ & $x_0^2$  & $6$  & $3.34\cdot 10^{-4067}$ & $6.45\cdot 10^{-945}$  & $1.45\cdot 10^{-38}$ & $2.51\cdot 10^{-1653}$
                           & $6.29\cdot 10^{-12}$   & $3.89\cdot 10^{-709}$  & $1.58\cdot 10^{-28}$\\ [0.6em]
 $(\star)$ & $x_0^3$  & $6$  & $2.32\cdot 10^{-6750}$ & $5.88\cdot 10^{-1688}$ & $3.27\cdot 10^{-47}$ & $1.20\cdot 10^{-2953}$
                           & $2.81\cdot 10^{-14}$   & $2.04\cdot 10^{-1266}$ & $3.27\cdot 10^{-47}$\\ [0.2em]
\hline\\
   FDN   & $x_0^1$  & $10$ & $1.83\cdot 10^{-5372}$ & $2.75\cdot 10^{-1791}$ & $8.88\cdot 10^{-13}$ & $3.94\cdot 10^{-2985}$
                           & $1.29\cdot 10^{-19}$   & $1.43\cdot 10^{-1194}$ & $1.72\cdot 10^{-34}$ \\ [0.6em]
$\rho=3$ & $x_0^2$  & $8$  & $9.36\cdot 10^{-6199}$ & $1.02\cdot 10^{-2066}$ & $3.78\cdot 10^{-14}$ & $3.50\cdot 10^{-3444}$
                           & $6.42\cdot 10^{-20}$   & $3.43\cdot 10^{-1378}$ & $1.02\cdot 10^{-34}$\\ [0.6em]
         & $x_0^3$  & $7$  & $5.20\cdot 10^{-3475}$ & $8.39\cdot 10^{-1159}$ & $2.28\cdot 10^{-13}$ & $5.44\cdot 10^{-1931}$
                           & $2.94\cdot 10^{-17}$   & $6.48\cdot 10^{-773}$  & $1.15\cdot 10^{-34}$\\ [0.2em]
\hline\hline
\end{tabular}}
\end{table}

\vspace{1mm}
Note that in Table \ref{angtable2} the method HMN is of $4$th order according to (\ref{FDN2}). Setting in ($\star$) $\rho=4$ the results obtained for ECOC and ACOC were excellent for the three initial conditions with $j=2$, but for COC it was necessary to take $j=20$. Hence in the computations shown in Table \ref{angtable2} was used $j=20$. In  ($\ast$), for NM method $j=3$ was used.

\subsubsection{Example 3}

A system of equations involving cubic polynomials,   $F_3(x)=0$, namely

\begin{center}
$\displaystyle \left\{ \begin{array}{r} x^3-3xy^2-1=0,\\ 3x^2y-y^3+1=0\,,\end{array}\right.$
\end{center}
is analyzed.

\vspace{1mm}
\noindent
\vspace{1mm}
Its solutions are  $(x,y)= {\left(5+3\sqrt{\mathstrut 3}\right)}^{1/3}\cdot\left(\dfrac{1}{2} ,\dfrac{\sqrt{\mathstrut 3}}{2}-1\right)$,
$(x,y)= {\left(5+3\sqrt{\mathstrut 3}\right)}^{1/3}\cdot\left(\dfrac{\sqrt{\mathstrut 3}}{2}-1, \dfrac{1}{2} \right)$
and $(x,y)= \left( -2^{-1/3},\, -2^{-1/3} \right)$.
We test the convergence of the methods towards the root
$\displaystyle \alpha = {\left(5+3\sqrt{\mathstrut 3}\right)}^{1/3}\cdot\left(\dfrac{\sqrt{\mathstrut 3}}{2}-1,\dfrac{1}{2}\right)^t
\approx(-0.2905145555, 1.084215081)^t\,.$

\vspace{2mm}
\noindent The initial values are
$x_0^1 =(-1,\,2)^t$    with $d_1=0.916$ and $D_1=10.0$\,; \
$x_0^2 =( -0.1,1.4)^t$  with $d_2=0.316$ and $D_2=1.702$\,; \ and
$x_0^3 =(-0.3,1.1)^t$ with $d_3=0.0158$ and $D_3=0.062$\,. Notice that in Table  \ref{angtable3} all the iterative methods give the same values of $\Delta{\overline{\rho}}_{k-1}$ and
$\Delta{\widehat{\rho}}_k\,$.

\begin{table}[ht]
\caption{Numerical results for system $F_3(x)=0$\,.}
\label{angtable3}
{\scriptsize\noindent
\begin{tabular}{cccccccccc}
         \hline\hline \\[-0.7em]
         &       &          &                      & \hfill COC & & \hfill ECOC & & \hfill ACOC & \\[0.2em]
 Method  & $x_0$   & $k$ & ${\|F_3(x_k)\|}_{\infty}$ & ${\|e_{k-1}\|}_{\infty}$             & $ \Delta{\overline{\rho}}_{k-1} $
                                                   & ${\|{\widetilde{e}}_{k}\|}_{\infty}$ & $ \Delta{\widetilde{\rho}}_k $
                                                   & ${\|\delta_{k}\|}_{\infty}$          & $ \Delta{\widehat{\rho}}_k $
                                                   \\[0.2em] \hline \\
    NM   & $x_0^1$  & $14$ & $9.85\cdot 10^{-3759}$ & $5.43\cdot 10^{-1880}$ & $9.00\cdot 10^{-5}$ & $2.82\cdot 10^{-2819}$
                           & $2.60\cdot 10^{-3}$    & $2.10\cdot 10^{-940}$  & $9.00\cdot 10^{-5}$ \\[0.6em]
$\rho=2$ & $x_0^2$  & $13$ & $3.52\cdot 10^{-4507}$ & $2.85\cdot 10^{-2244}$ & $3.68\cdot 10^{-4}$ & $4.57\cdot 10^{-3381}$
                           & $5.59\cdot 10^{-4}$    & $1.60\cdot 10^{-1127}$ & $3.68\cdot 10^{-4}$ \\[0.6em]
         & $x_0^3$  & $11$ & $3.97\cdot 10^{-3665}$ & $3.26\cdot 10^{-1833}$ & $6.29\cdot 10^{-6}$ & $1.77\cdot 10^{-2749}$
                           & $2.08\cdot 10^{-3}$    & $5.46\cdot 10^{-917}$  & $6.29\cdot 10^{-6}$ \\[0.2em]
\hline \\
   AMN   & $x_0^1$  & $9$  & $3.20\cdot 10^{-4152}$ & $9.37\cdot 10^{-1385}$ & $1.71\cdot 10^{-4}$ & $5.52\cdot 10^{-2307}$
                           & $3.88\cdot 10^{-3}$    & $1.97\cdot 10^{-923}$  & $1.71\cdot 10^{-4}$\\[0.6em]
$\rho=3$ & $x_0^2$  & $8$  & $8.73\cdot 10^{-3422}$ & $2.54\cdot 10^{-1141}$ & $4.38\cdot 10^{-4}$ & $1.89\cdot 10^{-1901}$
                           & $5.50\cdot 10^{-3}$    & $4.32\cdot 10^{-761}$  & $4.37\cdot 10^{-4}$\\[0.6em]
         & $x_0^3$  & $7$  & $2.23\cdot 10^{-3841}$ & $3.60\cdot 10^{-1281}$ & $1.03\cdot 10^{-3}$ & $4.08\cdot 10^{-2134}$
                           & $3.78\cdot 10^{-3}$    & $2.25\cdot 10^{-554}$  & $1.03\cdot 10^{-3}$\\[0.2em]
\hline\\
   HMN   & $x_0^1$  & $9$  & $5.18\cdot 10^{-7826}$ & $4.95\cdot 10^{-2609}$ & $4.43\cdot 10^{-4}$ & $5.76\cdot 10^{-4348}$
                           & $5.85\cdot 10^{-4}$    & $8.41\cdot 10^{-1740}$ & $4.43\cdot 10^{-4}$\\ [0.6em]
$\rho=3$ & $x_0^2$  & $8$  & $2.12\cdot 10^{-5383}$ & $7.58\cdot 10^{-1795}$ & $3.13\cdot 10^{-4}$ & $8.65\cdot 10^{-2991}$
                           & $1.27\cdot 10^{-3}$    & $4.55\cdot 10^{-1795}$ & $3.13\cdot 10^{-4}$\\ [0.6em]
 & $x_0^3\diamond$  & $7$  & $1.64\cdot 10^{-4740}$ & $1.20\cdot 10^{-1580}$ & $5.52\cdot 10^{-4}$ & $3.78\cdot 10^{-2633}$
                           & $4.56\cdot 10^{-4}$    & $2.47\cdot 10^{-1054}$ & $5.52\cdot 10^{-4}$\\ [0.2em]
\hline\\
   FDN   & $x_0^1$  & $8$  & $3.40\cdot 10^{-3218}$ & $1.58\cdot 10^{-1073}$ & $8.07\cdot 10^{-4}$ & $2.19\cdot 10^{-1788}$
                           & $4.39\cdot 10^{-3}$    & $8.23\cdot 10^{-716}$  & $8.07\cdot 10^{-4}$\\ [0.6em]
$\rho=3$ & $x_0^2$  & $8$  & $2.38\cdot 10^{-2869}$ & $2.87\cdot 10^{-957}$  & $9.83\cdot 10^{-4}$ & $8.18\cdot 10^{-1595}$
                           & $2.18\cdot 10^{-3}$    & $2.85\cdot 10^{-638}$  & $9.83\cdot 10^{-4}$\\ [0.6em]
         & $x_0^3$  & $7$  & $1.92\cdot 10^{-3591}$ & $5.43\cdot 10^{-1198}$ & $1.10\cdot 10^{-4}$ & $6.60\cdot 10^{-1996}$
                           & $2.76\cdot 10^{-3}$    & $8.45\cdot 10^{-799}$  & $1.10\cdot 10^{-4}$\\ [0.2em]
\hline\hline
\end{tabular}}
\end{table}

\subsubsection{Example 4}

We present an example consisting in the computation of the complex root of the equation $e^z=z$ with the smallest imaginary part. Here we solve the nonlinear system of equations

\begin{center}
$\displaystyle\left\{ \begin{array}{r} e^x\cos{y}=x,\\ e^x\sin{y}=y,\end{array}\right.\ $
\end{center}

\vspace{1mm}
\noindent that is obtained by considering the real and imaginary parts of the original equation $e^{x+iy}=x+iy$.

\vspace{1mm}
\noindent We apply the methods defined in previous sections setting
$x_0^1 =(0,2)^t$   with $d_1=0.6627$ and $D_1=1.091$\,; \
$x_0^2 =(0.2,1.1)^t$ with $d_2=0.2372$ and $D_2=0.354$\,; \ and
$x_0^3 =(0.3,1.3)^t$ with $d_3=0.03723$ and $D_3=0.0611$\,.

\begin{table}[ht]
\caption{Numerical results for system $F_4(x)=0$\,.}
\label{angtable4}
{\scriptsize \noindent
\begin{tabular}{cccccccccc}
         \hline\hline \\[-0.7em]
        &       &          &                      & \hfill COC & & \hfill ECOC & & \hfill ACOC & \\[0.2em]
 Method  & $x_0$   & $k$ & ${\|F_4(x_k)\|}_{\infty}$ & ${\|e_{k-1}\|}_{\infty}$             & $ \Delta{\overline{\rho}}_{k-1} $
                                                   & ${\|{\widetilde{e}}_{k}\|}_{\infty}$ & $ \Delta{\widetilde{\rho}}_k $
                                                   & ${\|\delta_{k}\|}_{\infty}$          & $ \Delta{\widehat{\rho}}_k $
                                                   \\[0.2em] \hline \\
    NM   & $x_0^1$ & $12$ & $8.06\cdot 10^{-4464}$ & $3.46\cdot 10^{-2232}$ & $1.54\cdot 10^{-5}$ & $4.92\cdot 10^{-3348}$
                          & $1.50\cdot 10^{-3}$    & $1.42\cdot 10^{-1116}$ & $1.54\cdot 10^{-5}$\\ [0.6em]
$\rho=2$ & $x_0^2$ & $12$ & $2.82\cdot 10^{-3616}$ & $2.00\cdot
10^{-1808}$ & $6.03\cdot 10^{-5}$ & $6.33\cdot 10^{-2712}$
                       & $9.55\cdot 10^{-4}(\ast)$ & $9.85\cdot 10^{-905}$  & $6.03\cdot 10^{-5}$\\ [0.6em]
         & $x_0^3$ & $11$ & $6.26\cdot 10^{-3517}$ & $9.42\cdot 10^{-1759}$ & $3.03\cdot 10^{-4}$ & $1.28\cdot 10^{-2637}$
                          & $9.10\cdot 10^{-4}$    & $7.31\cdot 10^{-880}$  & $3.03\cdot 10^{-4}$\\ [0.2em]
\hline \\
   AMN   & $x_0^1$ & $9$  & $1.58\cdot 10^{-3112}$ & $7.37\cdot 10^{-1038}$ & $3.28\cdot 10^{-7}$ & $1.11\cdot 10^{-1728}$
                          & $3.71\cdot 10^{-7}$    & $2.48\cdot 10^{-692}$  & $3.28\cdot 10^{-7}$\\ [0.6em]
$\rho=3$ & $x_0^2$ & $8$  & $1.89\cdot 10^{-5422}$ & $6.64\cdot
10^{-1808}$ & $9.64\cdot 10^{-5}$ & $7.56\cdot 10^{-3013}$
                          & $4.23\cdot 10^{-4}$    & $1.10\cdot 10^{-1205}$ & $9.64\cdot 10^{-5}$\\ [0.6em]
         & $x_0^3$ & $7$  & $1.30\cdot 10^{-3631}$ & $6.64\cdot 10^{-1211}$ & $2.16\cdot 10^{-4}$ & $3.59\cdot 10^{-2017}$
                      & $2.51\cdot 10^{-2}(\star)$ & $1.07\cdot 10^{-807}$  & $2.16\cdot 10^{-4}$\\ [0.2em]
\hline\\
   HMN   & $x_0^1$ & $9$  & $3.52\cdot 10^{-4442}$ & $5.29\cdot 10^{-1481}$ & $5.23\cdot 10^{-4}$ & $1.51\cdot 10^{-2466}$
                          & $6.13\cdot 10^{-3}$    & $5.38\cdot 10^{-988}$  & $5.23\cdot 10^{-4}$\\ [0.6em]
$\rho=3$ & $x_0^2$ & $7$  & $8.20\cdot 10^{-7487}$ & $8.96\cdot
10^{-2496}$ & $6.96\cdot 10^{-5}$ & $1.63\cdot 10^{-4159}$
                          & $1.33\cdot 10^{-3}$    & $1.82\cdot 10^{-1664}$ & $6.99\cdot 10^{-5}$\\ [0.6em]
         & $x_0^3$ & $7$  & $1.67\cdot 10^{-4248}$ & $1.90\cdot 10^{-2830}$ & $4.03\cdot 10^{-4}$ & $2.16\cdot 10^{-2360}$
                          & $4.37\cdot 10^{-3}$    & $6.83\cdot 10^{-2830}$ & $4.03\cdot 10^{-4}$\\ [0.2em]
\hline\\
   FDN   & $x_0^1$ & $8$  & $2.19\cdot 10^{-8244}$ & $1.60\cdot 10^{-2748}$ & $1.92\cdot 10^{-4}$ & $1.96\cdot 10^{-4580}$
                          & $1.20\cdot 10^{-3}$    & $1.23\cdot 10^{-1832}$ & $1.92\cdot 10^{-4}$\\ [0.6em]
$\rho=3$ & $x_0^2$ & $8$  & $9.95\cdot 10^{-4626}$ & $2.08\cdot
10^{-1542}$ & $1.40\cdot 10^{-3}$ & $3.35\cdot 10^{-2570}$
                          & $6.59\cdot 10^{-3}$    & $1.61\cdot 10^{-1028}$ & $1.40\cdot 10^{-4}$\\ [0.6em]
         & $x_0^3$ & $7$  & $4.67\cdot 10^{-3419}$ & $3.32\cdot 10^{-1140}$ & $1.39\cdot 10^{-3}$ & $1.46\cdot 10^{-633}$
                          & $8.44\cdot 10^{-4}$    & $1.05\cdot 10^{-1899}$ & $1.39\cdot 10^{-3}$\\ [0.2em]
\hline\hline
\end{tabular}}
\end{table}

\vspace{1mm}
\noindent In Table \ref{angtable4} in the cases $(\ast)$  and $(\star)$ was necessary to use  $j=3$  and $j=4$ respectively,  and again $\Delta{\overline{\rho}}_{k-1}$ and $\Delta{\widehat{\rho}}_k$ agree. This example can be found \cite{GGP}.

\subsubsection{Example 5}

A system $F_5(x)=0$ involving three nonlinear equations defined by

\begin{center}
$\displaystyle\left\{\begin{array}{r} xyz=1,\\ x+y-z^2=0,\\ x^2+y^2+z^2=9, \end{array}\right.\ $
\end{center}
is studied. We analyze the convergence of the methods towards the root
$$
\alpha \approx(2.14025812200,\, - 2.09029464225,\, - 0.22352512107)^t\,.
$$
The initial values are
$x_0^1 =(1.0,\,-1.0,\,0.1)^t$  with $d_1=1.14$ and $D_1=6.99$\,; \
$x_0^2 =(2.0,\,-2.0,\,0.0)^t$  with $d_2=0.224$ and  $D_3=1.00$\,; \  and  \
$x_0^3 =(2.1,\,-2.1,\,-0.2)^t$ with $d_3=0.0403$ and  $D_3=0.14$\,. This example can be found in \cite{HMT}.

\begin{table}[hb]
\caption{Numerical results for system
$F_5(x)=0$\,.}\label{angtable5} {\scriptsize\noindent
\begin{tabular}{cccccccccc}
         \hline\hline \\[-0.7em]
         &       &          &                      & \hfill COC & & \hfill ECOC & & \hfill ACOC & \\[0.2em]
 Method  & $x_0$   & $k$ & ${\|F_5(x_k)\|}_{\infty}$ & ${\|e_{k-1}\|}_{\infty}$             & $ \Delta{\overline{\rho}}_{k-1} $
                                                   & ${\|{\widetilde{e}}_{k}\|}_{\infty}$ & $ \Delta{\widetilde{\rho}}_k $
                                                   & ${\|\delta_{k}\|}_{\infty}$          & $ \Delta{\widehat{\rho}}_k $
                                                   \\[0.2em] \hline \\
    NM   & $x_0^1$ & $14$ & $5.21\cdot 10^{-2980}$ & $1.89\cdot 10^{-1490}$ & $1.01\cdot 10^{-4}$ & $1.81\cdot 10^{-2235}$
                          & $1.47\cdot 10^{-4}$    & $9.60\cdot 10^{-746}$  & $1.01\cdot 10^{-4}$\\[0.6em]
$\rho=2$ & $x_0^2$ & $12$ & $1.94\cdot 10^{-3939}$ & $2.17\cdot
10^{-1970}$ & $7.69\cdot 10^{-4}$ & $5.85\cdot 10^{-2955}$
                          & $1.23\cdot 10^{-3}$    & $5.21\cdot 10^{-986}$  & $7.69\cdot 10^{-4}$ \\[0.6em]
         & $x_0^3$ & $11$ & $7.33\cdot 10^{-3861}$ & $6.97\cdot 10^{-1931}$ & $1.57\cdot 10^{-3}$ & $5.28\cdot 10^{-2896}$
                          & $8.40\cdot 10^{-4}$ & $7.48\cdot 10^{-966}$ & $1.57\cdot 10^{-3}$ \\[0.2em]
\hline \\
   AMN   & $x_0^1$ & $9$  & $6.58\cdot 10^{-3481}$ & $1.02\cdot 10^{-1160}$ & $1.63\cdot 10^{-4}$ & $2.75\cdot 10^{-1934}$
                          & $3.42\cdot 10^{-4}$ & $2.71\cdot 10^{-774}$ & $1.63\cdot 10^{-4}$ \\[0.6em]
$\rho=3$ & $x_0^2$ & $8$  & $5.97\cdot 10^{-6440}$ & $4.99\cdot
10^{-2147}$ & $6.63\cdot 10^{-4}$ & $9.62\cdot 10^{-3578}$
                          & $1.27\cdot 10^{-3}$ & $7.16\cdot 10^{-1432}$ & $6.63\cdot 10^{-4}$\\[0.6em]
         & $x_0^3$ & $7$  & $2.70\cdot 10^{-4280}$ & $2.98\cdot 10^{-1427}$ & $1.67\cdot 10^{-4}$ & $2.27\cdot 10^{-2378}$
                          & $2.74\cdot 10^{-3}$ & $4.33\cdot 10^{-952}$ & $1.67\cdot 10^{-4}$\\[0.2em]
\hline\\
   HMN   & $x_0^1$ & $8$  & $5.99\cdot 10^{-2777}$ & $1.37\cdot 10^{-925}$  & $6.60\cdot 10^{-15}$ & $6.33\cdot 10^{-1543}$
                          & $9.51\cdot 10^{-13}$   & $4.63\cdot 10^{-618}$  & $5.28\cdot 10^{-31}$ \\ [0.6em]
$\rho=3$ & $x_0^2$ & $7$  & $1.67\cdot 10^{-3065}$ & $8.94\cdot
10^{-1022}$ & $4.27\cdot 10^{-15}$ & $3.11\cdot 10^{-1703}$
                          & $1.09\cdot 10^{-12}$   & $3.49\cdot 10^{-682}$  & $1.70\cdot 10^{-35}$\\ [0.6em]
         & $x_0^3$ & $7$  & $8.00\cdot 10^{-5214}$ & $6.99\cdot 10^{-1738}$ & $2.49\cdot 10^{-17}$ & $9.60\cdot 10^{-2897}$
                          & $1.15\cdot 10^{-13}$   & $1.37\cdot 10^{-1159}$ & $8.81\cdot 10^{-35}$\\ [0.2em]
\hline\\
   FDN   & $x_0^1$ & $16$ & $1.21\cdot 10^{-3001}$ & $4.46\cdot 10^{-1001}$ & $3.54\cdot 10^{-5}$ & $4.26\cdot 10^{-1668}$
                          & $7.99\cdot 10^{-5}$    & $9.55\cdot 10^{-668}$  & $3.54\cdot 10^{-5}$\\ [0.6em]
$\rho=3$ & $x_0^2$ & $8$  & $3.71\cdot 10^{-5272}$ & $6.10\cdot
10^{-1758}$ & $1.28\cdot 10^{-3}$ & $1.47\cdot 10^{-2929}$
                          & $9.54\cdot 10^{-4}$ & $1.85\cdot 10^{-1172}$ & $1.28\cdot 10^{-3}$ \\ [0.6em]
         & $x_0^3$ & $7$  & $7.53\cdot 10^{-3880}$ & $9.88\cdot 10^{-1294}$ & $1.45\cdot 10^{-4}$ & $8.35\cdot 10^{-2156}$
                          & $1.26\cdot 10^{-3}$ & $7.40\cdot 10^{-863}$ & $1.45\cdot 10^{-4}$\\ [0.2em]
\hline\hline
\end{tabular}}
\end{table}

\subsubsection{Example 6}

A system of four equations $F_6(x)=0$ found in \cite{CT}, and given by

$$\displaystyle\left\{\begin{array}{r} yz+t(y+z)=0,\\
xz+t(x+z)=0,\\ xy+t(x+y)=0,\\ xy+xz+yz=1,\end{array}\right.\ $$
is solved. Two roots of $F_6(x)=0$ are
$(x,y,z,t)=\pm\left(\dfrac{\sqrt{\mathstrut 3}}{3},\dfrac{\sqrt{\mathstrut 3}}{3},
\dfrac{\sqrt{\mathstrut 3}}{3},
\dfrac{-\sqrt{\mathstrut 3}}{6}\right)$. Here, we study the convergence of the methods towards
$$\displaystyle \alpha =\left(\dfrac{\sqrt{\mathstrut 3}}{3},\dfrac{\sqrt{\mathstrut 3}}{3},\dfrac{\sqrt{\mathstrut 3}}{3},
\dfrac{-\sqrt{\mathstrut 3}}{6}\right)^t \approx(0.5773502692, 0.5773502692, 0.5773502692,- 0.2886751346)^t\,.$$

\vspace{1mm}
\noindent Initial conditions are
$x_0^1 =(0.5,\,0.5,\,0.5,\,0.2)^t$  with  $d_1=0.4887$ and $D_1=0.45$; \
$x_0^2 =(0.55,\,0.55,\,0.55,\,-0.1)^t$  with $d_2=0.1887$ and $D_2=0.1925$;
 \  and \
$x_0^3 =(0.6,\,0.6,\,0.6,\,-0.3)^t$ with $d_3=0.02265$ and $D_3=0.08$. Like in the second example HMN method is of $4$th order $(\star)$. In Table \ref{angtable6} in the methods marked  $(\ast)$ $j=5$ was used.

\begin{table}[ht]
\caption{Numerical results for system $F_6(x)=0$\,.}
\label{angtable6} {\scriptsize \noindent
\begin{tabular}{cccccccccc}
         \hline\hline \\[-0.7em]
         &       &          &                      & \hfill COC & & \hfill ECOC & & \hfill ACOC & \\[0.2em]
 Method  & $x_0$   & $k$ & ${\|F_6(x_k)\|}_{\infty}$ & ${\|e_{k-1}\|}_{\infty}$             & $ \Delta{\overline{\rho}}_{k-1} $
                                                   & ${\|{\widetilde{e}}_{k}\|}_{\infty}$ & $ \Delta{\widetilde{\rho}}_k $
                                                   & ${\|\delta_{k}\|}_{\infty}$          & $ \Delta{\widehat{\rho}}_k $
                                                   \\[0.2em] \hline \\
    NM   & $x_0^1$ & $12$ & $3.20\cdot 10^{-4679}$ & $2.96\cdot 10^{-2339}$ & $5.14\cdot 10^{-4}$ & $2.28\cdot 10^{-3509}$
                          & $7.21\cdot 10^{-4}$  & $8.96\cdot 10^{-1172}$ & $5.14\cdot 10^{-4}$\\[0.6em]
$\rho=2$ & $x_0^2$ & $12$ & $3.75\cdot 10^{-3304}$ & $1.07\cdot
10^{-1650}$ & $7.28\cdot 10^{-4}$ & $2.52\cdot 10^{-2477}$
                          & $4.86\cdot 10^{-4}$  & $2.35\cdot 10^{-827}$ & $7.28\cdot 10^{-4}$\\[0.6em]
$(\ast)$ & $x_0^3$ & $11$ & $3.60\cdot 10^{-3514}$ & $1.10\cdot
10^{-1757}$ & $1.24\cdot 10^{-15}$ & $3.37\cdot 10^{-2636}$
                          & $8.96\cdot 10^{-18}$ & $3.08\cdot 10^{-879}$ & $4.72\cdot 10^{-33}$\\[0.2em]
\hline \\
     AMN & $x_0^1$ & $8$  & $4.38\cdot 10^{-7501}$ & $3.15\cdot 10^{-2498}$ & $1.72\cdot 10^{-3}$
                          & $1.67\cdot 10^{-4166}$  & $1.03\cdot 10^{-3}$ & $4.46\cdot 10^{-1668}$ & $1.72\cdot 10^{-3}$\\[0.6em]
$\rho=3$ & $x_0^2$ & $7$  & $1.16\cdot 10^{-3528}$ & $2.17\cdot
10^{-1174}$ & $3.65\cdot 10^{-3}$ & $7.51\cdot 10^{-1959}$
                          & $2.19\cdot 10^{-3} $ & $3.47\cdot 10^{-785}$ & $3.65\cdot 10^{-3}$\\[0.6em]
$(\ast)$ & $x_0^3$ & $7$  & $1.13\cdot 10^{-3752}$ & $1.63\cdot
10^{-1251}$ & $1.17\cdot 10^{-13}$ & $2.06\cdot 10^{-2085}$
                          & $1.32\cdot 10^{-21}$ & $1.26\cdot 10^{-834}$ & $4.71\cdot 10^{-35}$\\[0.2em]
\hline\\
       HMN  & $x_0^1$ & $6$ & $1.37\cdot 10^{-4681}$ & $3.30\cdot 10^{-1168}$ & $8.24\cdot 10^{-3}$ & $4.07\cdot 10^{-2046}$
                            & $4.71\cdot 10^{-3}$ & $1.23\cdot 10^{-878}$ & $8.24\cdot 10^{-3}$\\ [0.6em]
   $\rho=4$ & $x_0^2$ & $6$ & $1.30\cdot 10^{-6611}$ & $1.07\cdot 10^{-1650}$ & $5.83\cdot 10^{-3}$ & $1.73\cdot 10^{-2890}$
                            & $3.34\cdot 10^{-3}$    & $1.61\cdot 10^{-1240}$ & $5.83\cdot 10^{-3}$\\ [0.6em]
 $(\star)$   & $x_0^3\star$ & $6$ & $3.24\cdot 10^{-7028}$ & $1.10\cdot 10^{-1757}$ & $2.53\cdot 10^{-20}$ & $1.87\cdot 10^{-3075}$
                            & $2.40\cdot 10^{-25}$   & $1.71\cdot 10^{-1318}$ & $9.09\cdot 10^{-46}$
\\ [0.2em] \hline\\
   FDN   & $x_0^1$ & $8$ & $8.73\cdot 10^{-6401}$ & $1.50\cdot 10^{-2131}$ & $2.01\cdot 10^{-3}$& $2.48\cdot 10^{-3554}$
                         & $1.21\cdot 10^{-3}$    & $1.65\cdot 10^{-1423}$ & $2.01\cdot 10^{-3}$\\ [0.6em]
$\rho=3$ & $x_0^2$ & $7$ & $9.61\cdot 10^{-3188}$ & $7.57\cdot
10^{-1061}$ & $4.04\cdot 10^{-3}$ & $1.62\cdot 10^{-1769}$
                         & $2.43\cdot 10^{-3}$    & $2.15\cdot 10^{-709\ }$ & $4.04\cdot 10^{-3}$\\ [0.6em]
         & $x_0^3$ & $7$ & $1.43\cdot 10^{-3432}$ & $6.51\cdot 10^{-1145}$ & $3.72\cdot 10^{-16}$ & $1.20\cdot 10^{-1907}$
                         & $4.73\cdot 10^{-19}$   & $1.85\cdot 10^{-763\ }$ & $1.47\cdot 10^{-35}$\\ [0.2em]
\hline\hline
\end{tabular}}
\end{table}


\subsubsection{Example 7}

In \cite{GGP} the following boundary value problem
$$ y^{\prime\prime}+y^3=0\,,\quad y(0)=0\,,\quad y(1)=1\,,$$
was posed. To solve it, we consider the following partition of the interval $[0,1]:$
$$ u_0=0<u_1<u_2<\dots <u_{n-1}<u_n=1\,,\quad u_{j+1}=u_j+h\,,\ h=1/n\,.$$

Let us define $y_0=y(u_0)=0\,,\ y_1=y(u_1)\,,\dots\,, y_{n-1}=y(u_{n-1})\,,\ y_n=y(u_n)=1\,.$ If we discretize the problem by using the following numerical formula for the second derivative
$$y^{\prime\prime}_k \approx \, \dfrac{y_{k-1}-2y_k+2y_{k+1}}{h^2}\,,\quad k=2,3,\dots\,,\, n-1\,,$$

we obtain a $(n-1)\times(n-1)$ system of nonlinear equations:
$$\displaystyle\left\{\begin{array}{rr} 2y_1-h^2y^3_1-y_2=0\,,\\
                           -y_{k-1}+2y_k-h^2y^3_k-y_{k+1}=0\,, & k=2,3,\dots\,,\, n-1\,.\\
                            y_{n-2}+2y_{n-1}-h^2y^3_{n-1}=1\,,\end{array}\right.\ $$
In particular, we solve this problem for $n=10$. The solution $\alpha$ is

 $\displaystyle (.105541119905921, .211070483662496, .316505813937525, .421624081569127, .525992841283953,\\
 .628906344657317, .729332377591977, .825878904047790, .916792309006097)^t\,.$

\vspace{1mm}
\noindent We start with the points $x_0^1, x_0^2, x_0^3$, defined by
$x_0^1 =(1,0,-1,0,1,0,-1,0,1)^t$ with $d_1=1.7290$ and $D_1=1.990$;
$x_0^2 =(0, 0, 0, .5, .5, .5, 1, 1, 1)^t$  with $d_2=0.3165$  and $D_2=0.5012$ and
$x_0^3 =(0.1,0.2,0.3,0.4,0.5,0.6,0.7,0.8,0.9)^t$  with $d_3=0.02933$ and $D_3=0.007290$. In Table \ref{angtable7} we used $j=5$  in MN method $(\ast)$.

\begin{table}[ht]
\caption{Numerical results for system $F_7(x)=0$\,.}
\label{angtable7} {\scriptsize \noindent
\begin{tabular}{cccccccccc}
         \hline\hline \\[-0.7em]
         &       &          &                      & \hfill COC & & \hfill ECOC & & \hfill ACOC & \\[0.2em]
 Method  & $x_0$   & $k$ & ${\|F_7(x_k)\|}_{\infty}$ & ${\|e_{k-1}\|}_{\infty}$             & $ \Delta{\overline{\rho}}_{k-1} $
                                                   & ${\|{\widetilde{e}}_{k}\|}_{\infty}$ & $ \Delta{\widetilde{\rho}}_k $
                                                   & ${\|\delta_{k}\|}_{\infty}$          & $ \Delta{\widehat{\rho}}_k $
                                                    \\[0.2em] \hline \\
    NM   & $x_0^1$ & $12$ & $6.95\cdot 10^{-3714}$ & $1.92\cdot 10^{-1856}$ & $8.19\cdot 10^{-11}$ & $1.06\cdot 10^{-2784}$
                          & $1.02\cdot 10^{-10}$   & $5.52\cdot 10^{-929\ }$& $8.19\cdot 10^{-11}$ \\[0.6em]
$\rho=2$ & $x_0^2$ & $11$ & $7.42\cdot 10^{-3010}$ & $1.98\cdot
10^{-1504}$ & $4.50\cdot10^{-9}$   & $1.11\cdot 10^{-2256}$
                          & $1.81\cdot 10^{-9}$    & $5.61\cdot 10^{-753\ }$& $4.51\cdot 10^{-9}$  \\[0.6em]
$ (\ast)$  & $x_0^3$ & $11$ & $1.39\cdot 10^{-4759}$ & $2.72\cdot
10^{-2379}$ & $1.23\cdot 10^{-9}$  & $5.65\cdot 10^{-3569}$
                          & $1.51\cdot 10^{-9}$    & $2.08\cdot 10^{-1190}$ & $1.23\cdot 10^{-9}$  \\[0.2em]
\hline \\
   AMN   & $x_0^1$ & $8$  & $3.81\cdot10^{-4298}$ & $2.22\cdot10^{-1432}$ & $5.53\cdot10^{-6}$ & $2.24\cdot10^{-2387}$
                          & $1.91\cdot10^{-6}$    & $1.01\cdot10^{-955\ }$ & $5.53\cdot10^{-6}$ \\[0.6em]
$\rho=3$ & $x_0^2$ & $7$  & $2.35\cdot10^{-3175}$ &
$4.05\cdot10^{-1058}$ & $1.55\cdot10^{-5}$ & $1.32\cdot10^{-1763}$
                          & $1.18\cdot10^{-5}$    & $3.26\cdot10^{-706\ }$ & $1.55\cdot10^{-5}$ \\[0.6em]
         & $x_0^3$ & $7$  & $1.90\cdot10^{-5349}$ & $8.17\cdot10^{-1783}$ & $1.86\cdot10^{-5}$ & $1.98\cdot10^{-2971}$
                          & $2.57\cdot10^{-6}$    & $2.42\cdot10^{-1189}$ & $1.86\cdot10^{-5}$ \\[0.2em]
\hline\\
   HMN   & $x_0^1$ & $8$  & $1.23\cdot10^{-3764}$ & $1.35\cdot10^{-1254}$ & $2.68\cdot10^{-6}$ & $5.68\cdot10^{-2091}$
                          & $1.56\cdot10^{-6}$    & $4.20\cdot10^{-837}$  & $2.68\cdot10^{-6}$ \\[0.6em]
$\rho=3$ & $x_0^2$ & $7$  & $2.12\cdot10^{-3201}$ &
$7.52\cdot10^{-1067}$ & $3.78\cdot10^{-5}$ & $4.61\cdot10^{-1778}$
                          & $4.38\cdot10^{-5}$    & $6.13\cdot10^{-712}$  & $3.78\cdot10^{-5}$   \\[0.6em]
         & $x_0^3$ & $7$  & $3.45\cdot10^{-4877}$ & $1.90\cdot10^{-1625}$ & $3.45\cdot10^{-6}$ & $4.67\cdot10^{-2709}$
                          & $1.79\cdot10^{-6}$    & $2.45\cdot10^{-1084}$ & $3.45\cdot10^{-6}$   \\[0.2em]
\hline\\
   FDN   & $x_0^1$ & $8$  & $9.33\cdot10^{-6553}$ & $5.38\cdot10^{-2184}$ & $1.20\cdot10^{-7}$ & $6.10\cdot10^{-3640}$
                          & $2.81\cdot10^{-7}$    & $1.13\cdot10^{-1456}$ & $1.20\cdot10^{-7}$ \\[0.6em]
$\rho=3$ & $x_0^2$ & $8$  & $3.54\cdot10^{-3092}$ &
$1.81\cdot10^{-1030}$ & $1.47\cdot10^{-7}$ & $2.14\cdot10^{-1717}$
                          & $1.71\cdot10^{-6}$    & $1.18\cdot10^{-687\ }$& $1.46\cdot10^{-7}$\\[0.6em]
         & $x_0^3$ & $7$  & $6.47\cdot10^{-4760}$ & $2.21\cdot10^{-1586}$ & $3.62\cdot10^{-8}$ & $6.43\cdot10^{-2644}$
                          & $2.69\cdot10^{-8}$    & $2.91\cdot10^{-1058}$ & $3.62\cdot10^{-8}$\\[0.2em]
\hline\hline\\
\end{tabular}}
\end{table}

\vspace{5mm}
Finally, after performing the numerical computation of the previous examples, we have realized that there are not significant changes in the components of the approximations of the error vectors when using the above defined norm.


\section{Concluding remarks}

A generalization to several variables of a technique used to compute analytically the error equation of iterative methods without memory for one variable is presented. The key idea is to use formal power series.

\vspace{2mm}
So far, using the definition of COC for one variable, an iterative method is applied to a set of functions and numerical verification of the value of local order claimed in the error equation is performed. In this paper we generalize this procedure to several variables. Furthermore, when the root is not known, as usually occurs in real problems, we have overcome this situation introducing the numerical computation of the order ECOC and ACOC.

\vspace{2mm}
To obtain an approximation of the root with high precision it is necessary to define an adaptive float arithmetics allowing us, in each step, to get the appropriate length of the mantissa. For several variables this increasing value of the mantissa is defined according to the knowledge or not of the root.

\vspace{2mm}
To illustrate the technique presented  seven examples are worked and completely solved. In each one, four iterative methods have been carried out. Explicitly, agreements and differences are pointed out.

\vspace{2mm}
After comparing COC, ECOC and ACOC we think that ACOC is the best when the root is known or not due to its speed of convergence and its adaptability to the definition of the mantissa length.

%
%
%
%

\end{document}